\documentstyle[leqno]{book}
\newfont{\cirilrm}{wncyr10 scaled 1000}
\newfont{\cirilbf}{wncyb10 scaled 1000}
\newfont{\cirilsf}{wncyss10 scaled 1000}
\newfont{\cirilit}{wncyi10 scaled 1000}
\newfont{\cirilsc}{wncysc10 scaled 1000}

\newcounter{supersection}[section]
\newtheorem{th}[supersection]{Theorem}
\newtheorem{lm}[supersection]{Lemma}
\newtheorem{re}[supersection]{Remark}

\def\bibname{\textbf{REFERENCES}}
\def\thebibliography#1{\paragraph*{\uppercase{\bibname}}\list
{\bf [\arabic{enumi}]}{\settowidth\labelwidth{[#1]}\leftmargin\labelwidth
\advance\leftmargin\labelsep\usecounter{enumi}}
\def\newblock{\hskip .11em plus .33em minus .07em}
\sloppy\clubpenalty4000\widowpenalty4000
\sfcode`\.=1000\relax}

\def\F11{{}_{1}\mbox{\rm F}{}_{1}}

\def\stop{\mbox{\footnotesize {\vrule width 6pt height 6pt}}}

\begin{document}

\thispagestyle{plain}

\centerline{\large \bf SOME COMBINATORIAL ASPECTS}

\medskip

\centerline{\large \bf OF COMPOSITION OF A SET OF FUNCTIONS}
\footnotetext{Research partially supported by the MNTRS, Serbia \& Montenegro,
Grant No. 1861.}

\medskip

\begin{center}

\large \it Branko J. Male\v sevi\' c

\end{center}

\medskip

\begin{center}
\parbox{25.0cc}{\scriptsize \bf
Abstract. In this paper we determine a number of meaningful compositions
of higher order of a set of functions, which is considered in
{\bf \cite{Malesevic_98}}, in implicit and explicit form.
Results which are obtained are applied to the vector analysis
in order to determine  the number of meaningful differential operations
of higher order.}
\end{center}

\noindent
\section{\large \bf \boldmath \hspace*{-7.0 mm}
1. The composition of a set of functions}

Main topic of consideration in this paper is the set of  functions
$\mbox{$\cal A$}_{n}$, for \mbox{$n\!=\!2,3,\ldots\;$}, determined
in the following form:
\begin{equation}
\label{A_12}
\begin{tabular}{cc}
$\!\!\!\begin{array}{ll}
\mbox{\small $\mbox{$\cal A$}_{n}\;(n\!=\!2m)$:}\!\!
           & \mbox{\small $\nabla_{1}$}   : \mbox{A}_{0} \!\rightarrow\! \mbox{A}_{1}   \\
           & \mbox{\small $\nabla_{2}$}   : \mbox{A}_{1} \!\rightarrow\! \mbox{A}_{2}   \\
           & \,\,\vdots                                                                 \\
           & \mbox{\small $\nabla_{i}$}   : \mbox{A}_{i} \!\rightarrow\! \mbox{A}_{i+1} \\
           & \,\,\vdots                                                                 \\
           & \mbox{\small $\nabla_{m}$}   : \mbox{A}_{m-1} \!\rightarrow\! \mbox{A}_{m} \\
           & \mbox{\small $\nabla_{m+1}$} : \mbox{A}_{m} \!\rightarrow\! \mbox{A}_{m-1} \\
           & \,\,\vdots                                                                 \\
           & \mbox{\small $\nabla_{n-j}$} : \mbox{A}_{j+1} \!\rightarrow\! \mbox{A}_{j} \\
           & \,\,\vdots                                                                 \\
           & \mbox{\small $\nabla_{n-1}$} : \mbox{A}_{2} \!\rightarrow\! \mbox{A}_{1}   \\
           & \mbox{\small $\nabla_{n}$}   : \mbox{A}_{1} \!\rightarrow\! \mbox{A}_{0}
             \mbox{\normalsize ,}
\end{array}$
 $\!\!\!\!$&$\!\!\!\!$
$ \begin{array}{ll}
\mbox{\small $\mbox{$\cal A$}_{n}\;(n\!=\!2m\!+\!1)$:}\!\!
           & \mbox{\small $\nabla_{1}$}   : \mbox{A}_{0} \!\rightarrow\! \mbox{A}_{1}   \\
           & \mbox{\small $\nabla_{2}$}   : \mbox{A}_{1} \!\rightarrow\! \mbox{A}_{2}   \\
           & \,\,\vdots                                                                 \\
           & \mbox{\small $\nabla_{i}$}   : \mbox{A}_{i} \!\rightarrow\! \mbox{A}_{i+1} \\
           & \,\,\vdots                                                                 \\
           & \mbox{\small $\nabla_{m}$}   : \mbox{A}_{m-1} \!\rightarrow\! \mbox{A}_{m} \\
           & \mbox{\small $\nabla_{m+1}$} : \mbox{A}_{m} \!\rightarrow\! \mbox{A}_{m}   \\
           & \mbox{\small $\nabla_{m+2}$} : \mbox{A}_{m} \!\rightarrow\! \mbox{A}_{m-1} \\
           & \,\,\vdots                                                                 \\
           & \mbox{\small $\nabla_{n-j}$} : \mbox{A}_{j+1} \!\rightarrow\! \mbox{A}_{j} \\
           & \,\,\vdots                                                                 \\
           & \mbox{\small $\nabla_{n-1}$} : \mbox{A}_{2} \!\rightarrow\! \mbox{A}_{1}   \\
           & \mbox{\small $\nabla_{n}$}   : \mbox{A}_{1} \!\rightarrow\! \mbox{A}_{0}
             \mbox{\normalsize .}
\end{array}$
\end{tabular}
\end{equation}
Additionally, we make an assumption that $\mbox{A}_{i}$ are
non-empty sets, for $i\!=\!0,1,\ldots,m$, where  $m\!=\![n/2]$.
For each set of functions ${\cal A}_{n}$ we determine the number
of meaningful compositions of higher order in implicit and
explicit form. Let us define a binary relation $\rho$ "{\em to be
in composition}" with $\nabla_{i}\rho\nabla_{j}=1$ iff the
composition $\nabla_{j} \circ \nabla_{i}$ is meaningful for $i, j
\in \{1, 2, \ldots, n\}$. Let us form an adjacency matrix
$\mbox{\tt A} = [a_{ij}]$ of the graph, determined by relation
$\rho$, with
\begin{equation}
\label{A_ij}
a_{ij}
=
\left\{
\begin{array}{rcl}
1 & : & (j = i+1)    \vee   (i+j = n+1)                         \\[2.0 ex]
0 & : & (j \neq i+1) \wedge (i+j \neq n+1)
\end{array}
\right.
\end{equation}
for $i, j \in \{1, 2, \ldots, n\}$. Thus, on the basis of the article
{\bf \cite{Malesevic_98}}, some implicit formulas for the number of
meaningful compositions are given by the following statement.

\break

\begin{th}
\label{Th_1_1}
Let $P_{n}(\lambda) \!=\! |\mbox{\tt A} - \lambda \mbox{\tt I}| \!=\!
\alpha_{0} \lambda^{n} + \alpha_{1} \lambda^{n-1} + \ldots + \alpha_{n}$
be the characteristic polynomial of the matrix $\mbox{\tt A} = [a_{ij}]$,
determined by {\rm (\ref{A_ij})}, and  $v_n = [ \, 1 \, \cdots \, 1 \, ]_{1 \times n}$.
If we denote by $f(k)$ the number of meaningful composition of
$k^{\it \footnotesize th}\!$-order of functions from ${\cal A}_{n}$,
then the following formulas are true$:$
\begin{equation}
\label{Th_1_1_Form_1}
f(k) = v_n \cdot \mbox{\tt A}^{k-1} \cdot v^{T}_n
\end{equation}
and
\begin{equation}
\label{Th_1_1_Form_2}
\alpha_{0} f(k) + \alpha_{1} f(k-1) + \ldots + \alpha_{n} f(k-n) = 0
\quad (k>n).
\end{equation}
\end{th}

\begin{re}
Generally, let a graph $G$, with vertices $\nu_{1},\ldots,\nu_{n}$, be determined by adjacency
matrix $\mbox{\tt A}$ and let $P_{n}(\lambda) \!=\! |\mbox{\tt A} - \lambda \mbox{\tt I}|
\!=\! \alpha_{0} \lambda^{n} + \alpha_{1} \lambda^{n-1} + \ldots + \alpha_{n}$ be
the characteristic polynomial of the matrix $\mbox{\tt A}$. If we denote with
$a_{ij}^{(k)}$ the number of $\nu_{i}, \nu_{j}$$\,-\,$walks of length $k$ in the graph $G$,
then for every choice of $\nu_{i}$ and $\nu_{j}$, the sequence $a_{ij}^{(k)}$ satisfies
the same recurrent relation $(\ref{Th_1_1_Form_2})$. The previous statement
is the first problem in the section $8.6$, of the supplementary problems page,
of~the~book~{\rm \cite{West_04}}.
\end{re}

\noindent
\section{\large \bf \boldmath \hspace*{-7.0 mm}
2. Some explicit formulas for the number of composition}

In this part we give some explicit formulas for the number of
meaningful compositions of functions from the set ${\cal A}_{n}$.
The following statements are true.

\begin{lm}
\label{Lm_2_1} The characteristic polynomial $P_{n}(\lambda)$ of
matrix $\mbox{\tt A} = [a_{ij}]$, determined by {\rm
(\ref{A_ij})}, fulfills the following recurrent relation
\begin{equation}
\label{Lm_2_1_Form_1}
P_{n}(\lambda) = \lambda^2 {\big (} P_{n-2}(\lambda) - P_{n-4}(\lambda) {\big )}.
\end{equation}
\end{lm}
{\bf Proof.} Expanding the determinant $P_{n}(\lambda) =
|\mbox{\tt A} - \lambda \mbox{\tt I}|$ by first column we have
\begin{equation}
\label{F_01}
P_{n}(\lambda)
=
-\lambda C_{n-1}(\lambda) + (-1)^{n+1}D_{n-1}(\lambda),
\end{equation}
where $C_{n-1}(\lambda)$ and $D_{n-1}(\lambda)$ are suitable
minors of the elements $a_{11}$ i $a_{n1}$ of the determinant
$P_{n}(\lambda)$. Continuing the expansion of the determinant
$C_{n-1}(\lambda)$ by the ending row we can conclude that
\begin{equation}
\label{F_02}
C_{n-1}(\lambda) = -\lambda P_{n-2}(\lambda).
\end{equation}
Further, let us remark that determinant $D_{n}(\lambda)$ has minor $P_{n-3}(\lambda)$
as follows
\begin{equation}
D_{n}(\lambda)
=\mbox{\footnotesize
$\left|
\mbox{\begin{tabular}{cccccccc}
\multicolumn{1}{r}{$1$}        & \multicolumn{1}{c}{$0$} & \multicolumn{1}{c}{$0$} & \multicolumn{1}{c}{$\cdots$}           & \multicolumn{1}{c}{$0$} & \multicolumn{1}{c}{$0$} & \multicolumn{1}{r}{$0$}        & \multicolumn{1}{r}{$1$}      \\
\multicolumn{1}{r}{$-\lambda$} & \multicolumn{1}{c}{$1$} & \multicolumn{1}{c}{$0$} & \multicolumn{1}{c}{$\cdots$}           & \multicolumn{1}{c}{$0$} & \multicolumn{1}{c}{$0$} & \multicolumn{1}{r}{$1$}        & \multicolumn{1}{r}{$0$}      \\ \cline{2-6}
\multicolumn{0}{r|}{$0$}       & \multicolumn{1}{c}{$ $} & \multicolumn{1}{c}{$ $} & \multicolumn{1}{c}{$ $}                & \multicolumn{1}{c}{$ $} & \multicolumn{1}{c}{$ $} & \multicolumn{1}{|r}{$0$}       & \multicolumn{1}{r}{$0$}      \\
\multicolumn{0}{r|}{$0$}       & \multicolumn{1}{c}{$ $} & \multicolumn{1}{c}{$ $} & \multicolumn{1}{c}{$ $}                & \multicolumn{1}{c}{$ $} & \multicolumn{1}{c}{$ $} & \multicolumn{1}{|r}{$0$}       & \multicolumn{1}{r}{$0$}      \\
\multicolumn{0}{r|}{$\vdots$}  & \multicolumn{1}{c}{$ $} & \multicolumn{1}{c}{$ $} & \multicolumn{1}{c}{$P_{n-3}(\lambda)$} & \multicolumn{1}{c}{$ $} & \multicolumn{1}{c}{$ $} & \multicolumn{1}{|r}{$\vdots$}  & \multicolumn{1}{r}{$\vdots$} \\
\multicolumn{0}{r|}{$0$}       & \multicolumn{1}{c}{$ $} & \multicolumn{1}{c}{$ $} & \multicolumn{1}{c}{$ $}                & \multicolumn{1}{c}{$ $} & \multicolumn{1}{c}{$ $} & \multicolumn{1}{|r}{$0$}       & \multicolumn{1}{r}{$0$}      \\
\multicolumn{0}{r|}{$0$}       & \multicolumn{1}{c}{$ $} & \multicolumn{1}{c}{$ $} & \multicolumn{1}{c}{$ $}                & \multicolumn{1}{c}{$ $} & \multicolumn{1}{c}{$ $} & \multicolumn{1}{|r}{$1$}       & \multicolumn{1}{r}{$0$}      \\ \cline{2-6}
\multicolumn{1}{r}{$1$}        & \multicolumn{1}{c}{$0$} & \multicolumn{1}{c}{$0$} & \multicolumn{1}{c}{$\cdots$}           & \multicolumn{1}{c}{$0$} & \multicolumn{1}{c}{$0$} & \multicolumn{1}{r}{$-\lambda$} & \multicolumn{1}{r}{$1$}
\end{tabular}}
\right|$}.
\end{equation}
If in the previous determinant we multiply the first row by $\!-1$
and add it to the ${n}^{\rm \footnotesize th}$-row and then, if
in the next step, we expand determinant by ending column, we can
conclude
\begin{equation}
\label{F_03} D_{n}(\lambda)
=
(-1)^{n-1} \lambda^2 P_{n-3}(\lambda).
\end{equation}

\break

\noindent
On the basis of expansion (\ref{F_01}) and formulas (\ref{F_02}),
(\ref{F_03}) it is true that
\begin{equation}
\quad P_{n}(\lambda) = \lambda^2 {\big (} P_{n-2}(\lambda) -
P_{n-4}(\lambda) {\big )}.\;\stop
\end{equation}

\begin{lm}
\label{Lm_2_2}
Characteristic polynomial $P_{n}(\lambda)$ of the matrix $\mbox{\tt A} = [a_{ij}]$,
determined by {\rm (\ref{A_ij})}, has the following explicit representation
\begin{equation}
\label{Lm_2_2_Form_1} \label{P_explicit} \;\; P_{n}(\lambda)
=
\left\{
\begin{array}{ccl}

\displaystyle\sum\limits_{k=1}^{[\frac{n+2}{4}]+1}{(-1)^{k-1}
{\:\mbox{\scriptsize $\displaystyle\frac{n}{2}\!-\!k\!+\!2$}\:
\choose
\:\mbox{\scriptsize $k\!-\!1$}\:}

\lambda^{n-2k+2}}
\!\!\!\!&\!\!:\!\!&\!\! n\!=\!2m,                                              \\[2.0 ex]
\!\!\!\displaystyle\sum\limits_{k=1}^{[\frac{n+2}{4}]+2}{\!\!\!\!(-1)^{k-1}\!{\Bigg (}

\!{\:\mbox{\scriptsize $\displaystyle\frac{n\!+\!3}{2}\!-\!k$}\:
\choose
\mbox{\scriptsize $k\!-\!1$}}

\!+\!
{\:\mbox{\scriptsize $\displaystyle\frac{n\!+\!3}{2}\!-\!k$}\:
\choose
\!\!\mbox{\scriptsize $k\!-\!2$}\:} \! \lambda \! {\Bigg )} \lambda^{n-2k+2}}

\!\!\!&\!\!:\!\!&\!\! n\!=\!2m\!+\!1.\!\!\!\!
\end{array}
\right.
\end{equation}

\end{lm}
{\bf Proof.} Let us determine a few initial characteristic polynomials in the following forms:
\begin{equation}
\begin{array}{l}
P_{2}(\lambda) = \lambda ^{2} - 1
= \displaystyle\sum\limits_{k=1}^{2}{(-1)^{k-1} {\:\mbox{\scriptsize $3\!-\!k$}\:
\choose \:\mbox{\scriptsize $k\!-\!1$}\:} \lambda^{4-2k}},             \\[1.5 ex]
P_{4}(\lambda) = \lambda ^{4} - 2\,\lambda ^{2}
=
\displaystyle\sum\limits_{k=1}^{2}{(-1)^{k-1} {\:\mbox{\scriptsize $4\!-\!k$}\:
\choose \:\mbox{\scriptsize $k\!-\!1$}\:} \lambda^{6-2k}};
\end{array}
\end{equation}
and
\begin{equation}
\begin{array}{l}
\;\;
P_{3}(\lambda)
\!=\!\lambda ^{3}\!-\!\lambda ^{2}\!-\!\lambda
\!=\!\displaystyle\sum\limits_{k=1}^{3}{(-1)^{k-1} {\Bigg (}\!{\:\mbox{\scriptsize $3\!-\!k$}\: \choose \:\mbox{\scriptsize $k\!-\!1$}\:}\lambda^{5-2k}
\!+\!
{\:\mbox{\scriptsize $3\!-\!k$}\: \choose \:\mbox{\scriptsize $k\!-\!2$}\:} \lambda^{6-2k}\!{\Bigg )}},                                 \\[1.5 ex]
\;\;
P_{5}(\lambda)
\!=\!\lambda ^{5}\!-\!\lambda ^{4}\!-\!2\,\lambda ^{3}\!+\!\lambda ^{2}
\!
=\!\displaystyle\sum\limits_{k=1}^{3}{(-1)^{k-1} {\Bigg (}\!{\:\mbox{\scriptsize $4\!-\!k$}\: \choose \:\mbox{\scriptsize $k\!-\!1$}\:}\lambda^{7-2k}
\!+\!{\:\mbox{\scriptsize $4\!-\!k$}\: \choose \:\mbox{\scriptsize $k\!-\!2$}\:} \lambda^{8-2k}\!{\Bigg )}}.
\end{array}
\end{equation}
Then the statement of this lemma follows by mathematical induction
on the basis of the recurrent relation (\ref{Lm_2_1_Form_1}).~\stop

\smallskip
\noindent
From theorem \ref{Th_1_1} and lemma \ref{Lm_2_2} the following statement follows.
\begin{th}

\label{Th_2_3} Let $\mbox{\tt A} = [a_{ij}]$ be the matrix
determined by {\rm (\ref{A_ij})}. Then the number of meaningful
composition of $k^{\it \footnotesize th}\!$-order of functions
from ${\cal A}_{n}$ fulfills the recurrent relation {\rm
(\ref{Th_1_1_Form_2})}, whereas $\alpha_{i}$ $(i=0,1,\ldots,n)$
are coefficients of the characteristic polynomial $P_{n}(\lambda)$
determined by {\rm (\ref{Lm_2_2_Form_1})}.
\end{th}

\noindent
Further, the following general statement is true.

\begin{lm}
\label{Lm_2_4}
Let \mbox{$\mbox{\tt A}^{k} = [a_{ij}^{(k)}]$} be the
$k^{\it \footnotesize th}\!$-power of the matrix $\mbox{\tt A}
\!=\! [a_{i,j}] \!\in\! \mbox{\bf C}^{n \times n}$ $(k \!\in\! N)$
and~let
\begin{equation}
\label{Char_A}
P_{n}(\lambda) \!=\! |\mbox{\tt A} - \lambda \mbox{\tt I}| \!=\!
\alpha_{0} \lambda^{n} + \alpha_{1} \lambda^{n-1} + \ldots +
\alpha_{n},
\end{equation}
is characteristic polynomial $P_{n}(\lambda)$ of the matrix
$\mbox{\tt A}$. If for each pair of indexes
$(i,j)\!\in\!\{1,2,\ldots,n\}^{2}$ the sequence $g_{ij}(m)$, for
$m > n$, is determined as a solution of the recurrent relation
\begin{equation}
\label{Rec_ij}
\alpha_{0} g_{ij}(m) + \alpha_{1} g_{ij}(m-1) + \ldots + \alpha_{n} g_{ij}(m-n)
=
0,
\end{equation}
on the basis of initial values $g_{ij}(1) = a_{ij}^{(1)},\,
g_{ij}(2) = a_{ij}^{(2)}, \ldots, g_{ij}(n) = a_{ij}^{(n)}$, then
matrix $\mbox{\tt G}_{m} = [g_{ij}(m)] \in \mbox{\bf C}^{n \times
n}$ is the $m^{\it \footnotesize th}\!$-power of the matrix
\mbox{\tt A} $(m \!\in\! N)$.
\end{lm}

\break

\noindent
{\bf Proof.} We prove equality $G_{m} = A^{m}$ by total mathematical induction
over $m \!\in\! N$. Indeed, for $m = 1, \ldots, n$ statement is true. Let
$m > n$. Let us assume that $G_{k} = A^{k}$ is true for each $k < m$. Then
for $k = m$, let us note $g_{ij}(k)$ fulfils
\begin{equation}
g_{ij}(m)
=
-\displaystyle\frac{1}{\alpha_{0}} {\Big (}\alpha_{1}g_{ij}(m-1) +
\ldots + \alpha_{n} g_{ij}(m-n) {\Big )},
\end{equation}
where  $\alpha_{0}\!=\!(-1)^{n}\!$. From the previous equality, on
the basis of {\sc Cayle}-{\sc Hamilton}'s theorem, it follows that
\begin{equation}
\begin{array}{rcl}
G_{m}

&\!\!=\!\!&
-\displaystyle\frac{1}{\alpha_{0}}

{\Big (} \alpha_{1}G_{m-1} + \ldots + \alpha_{n} G_{m-n} {\Big )}              \\[2.0 ex]
&\!\!=\!\!&
-\displaystyle\frac{1}{\alpha_{0}}

{\Big (} \alpha_{1}A^{m-1} + \ldots + \alpha_{n} A^{m-n} {\Big )}
=
A^{m}.\,~\stop
\end{array}
\end{equation}

\medskip
\noindent
On the basis of theorem \ref{Th_1_1} and lemmas \ref{Lm_2_4}, \ref{Lm_2_2}
the following statement follows.

\begin{th}
\label{Th_2_5}
Let $\mbox{\tt A} = [a_{ij}]$ be the matrix determined by {\rm (\ref{A_ij})} and let
\mbox{$\mbox{\tt A}^{m} = [a_{ij}^{(m)}]$} is the $m^{\it \footnotesize th}\!$-power
of the matrix $\mbox{\tt A}$ determined for each pair of the indexes
$(i,j)\!\in\!\{1,2,\ldots,n\}^{2}$, for $m\!>\!n$, by an explicit form
of the elements $a_{ij}^{(m)}$ on the basis of a recurrent relation
\begin{equation}
\alpha_{0} a_{ij}^{(m)} + \alpha_{1} a_{ij}^{(m-1)} + \ldots +
\alpha_{n} a_{ij}^{(m-n)}
=
0.
\end{equation}
For all that $a_{ij}^{(1)},\, a_{ij}^{(2)}, \ldots, a_{ij}^{(n)}$ are initial values for
the previous recurrent relation and $\alpha_{i}$ $(i=0,1,\ldots,n)$ are coefficients
of the characteristic polynomial $P_{n}(\lambda)$ determined by {\rm (\ref{Lm_2_2_Form_1})}.
Then, by the formula {\rm (\ref{Th_1_1_Form_1})}, number $f(k)$ of the meaningful composition
of $k^{\it \footnotesize th}\!$-order of functions over ${\cal A}_{n}$, is explicitly
determined.
\end{th}

\noindent
\section{\large \bf \boldmath \hspace*{-7.0 mm}
3. Examples from vector analysis}

We present some examples of counting the numbers of meaningful
differential operations of higher order in vector analysis
according to {\bf \cite{Malesevic_98}}. Let us start with the sets
of functions
\begin{equation}
\mbox{\rm A}_{i} = \{ \mbox{f} : \mbox{\bf R}^{n} \longrightarrow
\mbox{\bf R}^{n \choose i} \, | \, f_{1}, \ldots, f_{n\choose i}
\in C^{\infty}(\mbox{\bf R}^{n})\},
\end{equation}
for $i = 0, 1, \ldots, m$, where $m=[n/2]$. Let $\Omega^r
(\mbox{\bf R}^n)$ be the space (module) of differential forms of
degree $r = 0, 1, \ldots, n$ on the space $\mbox{\bf R}^n$ over
ring $\mbox{\rm A}_{0}$. For each $r$ let us choose the order of
basis elements $d x_{i_1} \! \wedge \ldots \wedge d x_{i_r}$. For
each $i$ let us determine
\begin{equation}
\label{Isomorph} \quad \varphi_{i} : \Omega^{i}(\mbox{\bf R}^{n})
\rightarrow \mbox{\rm A}_{i} \, (\mbox{\small $0$} \! \leq \! i \!
\leq \! \mbox{\small $m$}) \;\;\; \mbox{and} \;\;\; \varphi_{n-i}
: \Omega^{n-i}(\mbox{\bf R}^{n}) \rightarrow \mbox{\rm A}_{i} \,
(\mbox{\small $0$} \! \leq \! i \! < \! \mbox{\small $n\!-\!m$}),
\end{equation}
by forming $n\choose i$-tuple of coefficients with respect to the
basis elements in the given order. Let us remark that, in
comparison to {\bf \cite{Malesevic_98}}, we additionally consider
the order of basis elements. Consequently, previously introduced
functions are isomorphisms from the space of differential forms
into the space of vector functions. Next, we use the well-known
fact that $\Omega^i(\mbox{\bf R}^n)$ and $\Omega^{n-i}(\mbox{\bf
R}^n)$ are spaces of the same dimension $n \choose i$, for $i = 0,
1, \ldots, m$. They can be identified with $\mbox{\rm A}_i$, using
corresponding isomorphism (\ref{Isomorph}). Let's define {\em
differential operations of the first order} via exterior
differentiation operator $d$ as follows
\begin{equation}
\nabla_{r} = \varphi_{r} \circ d \circ \varphi_{r-1}^{-1} \quad
(\mbox{\small $1$} \leq r \leq \mbox{\small $n$}).
\end{equation}
Thus the following diagrams commute:
\begin{equation}
\mbox{{\normalsize \unitlength 0.80 mm 
\quad \begin{picture}(0,10)(48,10)                                                  
\put(7.0,-3){\scriptsize $(1 \leq i \leq m)$}                                       
\put(0,16.2){$\Omega^{\mbox{\tiny $i\!\!-\!\!1$}}$}                                 
\put(1,4){\vector(0,1){11}}                                                         
\put(2,9){\scriptsize $\varphi^{-1}_{\mbox{\tiny $i\!\!-\!\!1$}}$}                  
\put(0,0.8){$A_{\mbox{\tiny $i\!-\!1$}}$}                                           
\put(3,3){\vector(1,0){26}}                                                         
\put(14.8,4){\scriptsize $\nabla_{\mbox{\tiny $i$}}$}                               
\put(29.0,0.8){$A_{\mbox{\tiny $i$}}$}                                              
\put(30,15){\vector(0,-1){11}}                                                      
\put(31,9){\scriptsize $\varphi_{\mbox{\tiny $i$}}$}                                
\put(29.2,16.2){$\Omega^{\mbox{\tiny $i$}}$}                                        
\put(3,17){\vector(1,0){26}}                                                        
\put(15.3,18){\scriptsize $d$}                                                      
\put(65.0,-3){\scriptsize $(0 \leq i < n\!-\!m)$}                                   
\put(59.5,16.2){$\Omega^{\mbox{\tiny $n\!\!-\!\!(\!i\!\!+\!\!1\!)$}}$}              
\put(61,4){\vector(0,1){11}}                                                        
\put(62,9){\scriptsize $\varphi^{-1}_{\mbox{\tiny $n\!\!-\!\!(\!i\!\!+\!\!1\!)$}}$} 
\put(60,0.8){$A_{\mbox{\tiny $i\!\!+\!\!1$}}$}                                      
\put(63,3){\vector(1,0){26}}                                                        
\put(73,4){\scriptsize $\nabla_{\mbox{\tiny $n\!\!-\!\!i$}}$}                       
\put(89.0,0.8){$A_{\mbox{\tiny $i$}}$}                                              
\put(90,15){\vector(0,-1){11}}                                                      
\put(91,9){\scriptsize $\varphi_{\mbox{\tiny $n\!\!-\!\!i$}}$}                      
\put(89.2,16.2){$\Omega^{\mbox{\tiny $n\!\!-\!\!i$}}$}                              
\put(63,17){\vector(1,0){26}}                                                       
\put(74,18){\scriptsize $d$}                                                        
\end{picture} }} 
\end{equation}

\vspace*{10.0 mm}

\noindent
Hence, the differential operations $\nabla_{r}$ determine functions so that (\ref{A_12})
is fulfilled. Let us define {\em differential operations of the higher order}
as meaningful compositions of higher order of functions from the set
${\cal A}_{n} = \{\nabla_{1}, \ldots, \nabla_{n}\}$. Let us consider,
in the next sections, concrete dimensions $n=3,4,5,\ldots,10$.

\bigskip
\noindent
{\bf Three-dimensional vector analysis.} In the real three-dimensional space $\mbox{\bf R}^3$
we consider the following sets
\begin{equation}
\mbox{\rm A}_{0} \!=\! \{f\!:\!\mbox{\bf R}^3 \!\longrightarrow\!
\mbox{\bf R} \, | \, f\!\in\!C^{\infty}(\mbox{\bf R}^{3}) \} \;
\mbox{and}\; \mbox{\rm A}_{1} \!=\! \{\vec{f}\!:\!\mbox{\bf R}^3
\!\longrightarrow\! \mbox{\bf R}^{3} \, | \,
\vec{f}\!\in\!\vec{C}^{\infty}(\mbox{\bf R}^{3})  \}.
\end{equation}
Let $dx$, $dy$, $dz$ respectively be the basis vectors of the
space of $1$-forms and let \mbox{$dy \wedge dz$}, \mbox{$dz \wedge
dx$}, \mbox{$dx \wedge dy$} respectively be the basis vectors of
the space $2$-forms. Thus over the sets $\mbox{\rm A}_{0}$ and
$\mbox{\rm A}_{1}$ there exist $m=3$ differential operations of
the first order
\begin{equation}
\begin{array}{l}
\nabla_1 \, f = \mbox{\rm grad} \, f
: \mbox{\rm A}_{0} \longrightarrow \mbox{\rm A}_{1},                           \\[2.0 ex]
\nabla_2 \, \vec{f} \!=\! \mbox{\rm curl} \, \vec{f}

: \mbox{\rm A}_{1} \longrightarrow \mbox{\rm A}_{1},                           \\[2.0 ex]
\nabla_3 \, \vec{f} = \mbox{\rm div} \, \vec{f}

: \mbox{\rm A}_{1} \longrightarrow \mbox{\rm A}_{0}.
\end{array}
\end{equation}
Under the previous choice of order of basis vectors of spaces of $1$-forms and $2$-forms,
the previously defined operations of the first order coincide with differential operations
of first order in the classical vector analysis. Next, as a well-known fact,
there are $m=5$ differential operations of the second order. In the article
{\bf \cite{Malesevic_96}} it is proved that there exists $m=8$ differential
operations of the third order. Further, in the article {\bf \cite{Malesevic_98}},
it is proved that there exist $F_{k+3}$ differential operations of the
${k}^{\mbox{\rm \footnotesize th}}\!$-order, where $F_{k}$ is {\sc Fibonacci}'s number
of order $k$. Here we give the proof of the previous statement, on the basic of results
from the second part of this paper. Namely, using theorems \ref{Th_2_3} and
\ref{Th_2_5}, matrix $\mbox{\tt A}^{k}$ has the following explicit form
\begin{equation}
\mbox{\tt A}^{k}
=
\left[
\begin{array}{lll}

F_{k-1} & F_{k}   & F_{k}  \\
F_{k-1} & F_{k}   & F_{k}  \\
F_{k-2} & F_{k-1} & F_{k-1}

\end{array}
\right].
\end{equation}
Hence, using (\ref{Th_1_1_Form_1}), the number of differential operation of
the ${k}^{\rm \footnotesize th}\!$-order is determined by
\begin{equation}
f(k)
=
v_{3} \cdot \mbox{\tt A}^{k-1} \cdot v_{3}^{T}
=
F_{k-3} + 4 F_{k-2} + 4F_{k-1}
=
F_{k+3}.
\end{equation}

\bigskip
\noindent
{\bf Multidimensional vector analysis.} In the real $n$-dimensional space
$\mbox{\bf R}^{n}$ the number of differential operations is determined by
corresponding recurrent formulas, which for dimension $n = 3, 4, 5, \ldots, 10$,
we cite according to {\bf \cite{Malesevic_98}}:

{\footnotesize
\begin{center}

\begin{tabular}{|c|c|}
\hline {\small \rm dimension}:
 &     {\small \rm \quad recurrent relations for the number of meaningful operations: \quad}

                                                                                   \\ \hline
$n = \;$ 3 & $f(i+2)=f(i+1)+f(i)$                                                  \\ \hline
$n = \;$ 4 & $f(i+2)=2 f(i)$                                                       \\ \hline
$n = \;$ 5 & $f(i+3)=f(i+2) + 2 f(i+1) - f(i)$                                     \\ \hline
$n = \;$ 6 & $f(i+4)=3 f(i+2) - f(i)$                                              \\ \hline
$n = \;$ 7 & $f(i+5)=f(i+3) + 3 f(i+2) - 2  f(i+1) - f(i)$                         \\ \hline
$n = \;$ 8 & $f(i+4)=4 f(i+2) - 3 f(i)$                                            \\ \hline
$n = \;$ 9 & $f(i+5)=f(i+4) + 4 f(i+3) - 3 f(i+2) - 3 f(i+1) + f(i)$               \\ \hline
$n = 10$ & $f(i+6)=5 f(i+4) - 6 f(i+2) + f(i)$                                     \\ \hline
\end{tabular}

\end{center}}

\noindent
For dimensions $n = 3$, as we have shown in the previous consideration, the numbers
of differential operations of higher order are determined via {\sc Fibonacci} numbers.
Also, this is true for dimension $n=6$. Namely, using theorems \ref{Th_2_3} and \ref{Th_2_5},
matrix $\mbox{\tt A}^{k}$ has following explicit form
\begin{equation}
\mbox{\tt A}^{k}
=
\left\{
\begin{array}{ccl}
\mbox{\small $\left[
\begin{array}{llllll}
F_{2p-1} \!&\! 0        \!&\! F_{2p}   \!&\! 0        \!&\! F_{2p}   \!&\! 0        \\[1.5 ex]
0        \!&\! F_{2p}   \!&\! 0        \!&\! F_{2p-1} \!&\! 0        \!&\! F_{2p}   \\[1.5 ex]
F_{2p-2} \!&\! 0        \!&\! F_{2p-1} \!&\! 0        \!&\! F_{2p-1} \!&\! 0        \\[1.5 ex]
0        \!&\! F_{2p}   \!&\! 0        \!&\! F_{2p-1} \!&\! 0        \!&\! F_{2p}   \\[1.5 ex]
F_{2p-1} \!&\! 0        \!&\! F_{2p}   \!&\! 0        \!&\! F_{2p}   \!&\! 0        \\[1.5 ex]
0        \!&\! F_{2p-1} \!&\! 0        \!&\! F_{2p-2} \!&\! 0        \!&\! F_{2p-1}
\end{array}
\right]$}
&\!\!\!\!:& k=2p,                                                                   \\[12.0 ex]
\mbox{$\left[
\begin{array}{llllll}
0        \!&\! F_{2p+1} \!&\! 0        \!&\! F_{2p}   \!&\! 0        \!&\! F_{2p+1} \\[1.5 ex]
F_{2p}   \!&\! 0        \!&\! F_{2p+1} \!&\! 0        \!&\! F_{2p+1} \!&\! 0        \\[1.5 ex]
0        \!&\! F_{2p}   \!&\! 0        \!&\! F_{2p-1} \!&\! 0        \!&\! F_{2p}   \\[1.5 ex]
F_{2p}   \!&\! 0        \!&\! F_{2p+1} \!&\! 0        \!&\! F_{2p+1} \!&\! 0        \\[1.5 ex]
0        \!&\! F_{2p+1} \!&\! 0        \!&\! F_{2p}   \!&\! 0        \!&\! F_{2p+1} \\[1.5 ex]
F_{2p-1} \!&\! 0        \!&\! F_{2p}   \!&\! 0        \!&\! F_{2p}   \!&\! 0
\end{array}
\right]$}
&\!\!\!\!:\!\!\!\!& k=2p+1.
\end{array}
\right.
\end{equation}
Hence, using (\ref{Th_1_1_Form_1}), the number of  differential
operations of the ${k}^{\rm \footnotesize th}\!$-order is
determined as follows
\begin{equation}
f(k) = 2 \cdot F_{k+3}.
\end{equation}
For other dimensions $n=4,5,7,8,9,10$ the roots of suitable characteristic polynomials
are not related to {\sc Fibonacci} numbers.

\smallskip \noindent
Finally, let us outline that for all dimensions $n=3,4,5,6,7,8,9,10$
the values of the function $f(k)$, for initial values of the argument~$k$,
are given in {\bf \cite{Sloane_03}} as sequences $A020701$ $(n=3)$,
$A090989$ $(n=4)$, $A090990$ $(n=5)$, $A090991$ $(n=6)$,
$A090992$ $(n=7)$, $A090993$ $(n=8)$, $A090994$ $(n=9)$,
$A090995$ $(n=10)$ respectively.

\break

\bigskip

{\small
\noindent University of Belgrade,
            \hfill (Received 31. 03. 2004.)          \break
\noindent Faculty of Electrical Engineering,   \hfill\break
\noindent P.O.Box 35-54, $11120$ Belgrade,     \hfill\break
\noindent Serbia \& Montenegro                 \hfill\break
\noindent {\footnotesize \bf malesevic@kiklop.etf.bg.ac.yu}
\hfill}

\end{document}